\documentclass{article}    
\usepackage[final]{graphicx}   
\usepackage{psfrag}   
\usepackage{amsmath,amsfonts,amssymb,amsxtra}

\newcommand {\Real}{\ensuremath{{\mathbb{R}}}}
\newcommand {\Natural}{\ensuremath{{\mathbb{N}}}}
\newcommand {\Complex}{\ensuremath{{\mathbb{C}}}}

\newcommand{\A}{\ensuremath{\mathcal A}}

\newcommand{\X}{\ensuremath{\mathcal X}}

\newcommand{\N}{\ensuremath{\mathcal N}}

\newcommand{\vi}{\ensuremath{{\mathbf{v}}}}
\newcommand{\ex}{\ensuremath{{\mathbf{x}}}}

\newcommand{\one}{\ensuremath{{\mathbf{1}}}}

\newcommand{\dabilyu}{\ensuremath{{\mathbf{w}}}}

\newtheorem{theorem}{Theorem}

\newtheorem{lemma}{Lemma}
\newtheorem{fact}{Fact}
\newtheorem{algorithm}{Algorithm}

\newenvironment{proof}{\noindent {\bf Proof.}}{\hfill \hspace*{1pt}\hfill$\blacksquare$}

\begin{document}

\title{Synchronizing discrete-time neutrally stable linear systems via partial-state coupling} 
\author{S. Emre Tuna\\
{\small {\tt tuna@eee.metu.edu.tr}} }
\maketitle


\begin{abstract}                          
A basic result in synchronization of linear systems via output
coupling is presented.  For identical discrete-time linear systems
that are detectable from their outputs and neutrally stable, it is
shown that a linear output feedback law exists under which the coupled
systems globally asymptotically synchronize for all fixed connected
(asymmetrical) network topologies. An algorithm is provided to compute
such feedback law based on individual system parameters. A dual problem
is also presented and solved.
\end{abstract}

\section{Introduction}

A notable meeting point for many researchers from different fields is
the topic {\em synchronization}. One of the reasons for that comes
from the nature as synchronization in large networks of dynamical
systems is a frequently encountered phenomenon in biology. Among many
others, one can count synchronously discharging neurons, crickets
chirping in accord, and metabolic synchrony in yeast cell
suspensions. Another reason is the abundance of technological
applications: coupled synchronized lasers, vehicle formations, and
sensor networks, just to name a few. We refer the reader to the
surveys \cite{wang02,strogatz01,boccaletti06,olfati07} for references
and more examples.

The main issue in studying the synchronization of coupled dynamical
systems is the stability of synchronization. As in all cases where
stability is the issue, the question whose answer is sought is {\em
Under what conditions} will the individual systems synchronize? In a
simplified yet widely-studied scenario, where the individual system
dynamics are identical and the coupling between them is linear,
studies focus on two ingredients: the dynamics of an individual system
and the network topology. Starting with the {\em agreement algorithm}
in \cite{tsitsiklis86} a number of contributions
\cite{jadbabaie03,moreau05,ren05,angeli06,blondel05,olfati04} have
gathered around the case where the weakest possible assumptions are
made on the network topology at the expense of restrictive individual
system dynamics. It was established in those works on {\em multi-agent
systems} that when the individual system is taken to be an integrator
and the coupling is of full-state, synchronization ({\em consensus})
results for time-varying interconnections whose unions\footnote{By
{\em union of interconnections} we actually mean the union of the
graphs representing the interconnections.} over an interval are
assumed to be connected instead of that each 
interconnection at every instant is connected.
         
Another school of research investigates networks with more complicated
(nonlinear) individual system dynamics. When that is the case, the
restrictions on the network topology have to be made stricter in order
to ensure stability of synchronization. Generally speaking, more than
mere connectedness of the network has been needed: coupling strength
is required to be larger than some threshold and sometimes a
symmetry\footnote{A network is {\em symmetric} if the matrix
representing it is symmetric.}  or balancedness assumption is made on
the connection graph. Different (though related) approaches have
provided different insights over the years. The primary of such
approaches is based on the calculations of the eigenvalues of the
connection matrix and a parameter (e.g. the maximal Lyapunov exponent)
depending on the individual system dynamics
\cite{wu95,pecora98}. In endeavor to better understand
synchronization stability, tools from systems theory such as Lyapunov
functions \cite{belykh06,hui07}, passivity
\cite{arcak07,stan07}, contraction theory
\cite{slotine04}, and incremental input to state stability
($\delta$ISS) theory \cite{cai06} have also proved useful.

This paper studies a broad class of linear systems under weak
assumptions on the coupling structure and generalizes some of the
existing results on synchronization. Namely, we consider identical
individual discrete-time linear systems interacting via (diffusive)
output coupling under a fixed (time-invariant) network topology. The
contribution of the paper is in proving (via construction) the
following basic result, which seems to have been missing from the
literature. For a linear system\footnote{$x^{+}=Ax+u;\ y=Cx.$} that is
neutrally stable and detectable from its output, there always exists a
linear output feedback law that ensures the global asymptotic
synchronization of any connected (not necessarily symmetric nor
balanced) network of any number of coupled replicas of that system. To
fortify our contribution practically, we provide an algorithm to
compute one such feedback law. Solving that problem will also yield us
a solution to the dual problem where the coupling is of full state but
the input injection to the system is through a (non-identity) $B$
matrix. As expected for the dual case, detectability assumption will
need to be replaced by stabilizability. It is worth noting that our
main theorem makes a compromise result between the two previously
mentioned cases (i) where synchronization is established for very
primitive individual system dynamics, such as that of an integrator,
but under the weakest conditions on the network topology and (ii)
where the network topology has to satisfy stronger conditions, such as
that the coupling strength should be above a threshold, for want of
achieving synchronization for nonlinear individual system dynamics.

The remainder of the paper is organized as follows. Notation and
definitions reside in the next section. We give the problem statement
along with our assumptions in Section~\ref{sec:problemstatement}. In
Section~\ref{sec:pre} we provide a preliminary synchronization result
on a network of linear systems with orthogonal system matrices.  Then
we generalize that result to establish our main theorem in
Section~\ref{sec:main}. We work out the dual case in
Section~\ref{sec:dual}.
  
\section{Notation and definitions}

The number of elements in a (finite) set ${\mathcal S}$ is denoted by
$\#{\mathcal S}$. Let $\Natural$ denote the set of nonnegative
integers. Let $|\cdot|$ denote 2-norm.  Identity matrix in
$\Real^{n\times n}$ is denoted by $I_{n}$. A matrix
$Q\in\Real^{n\times n}$ is {\em orthogonal} if
$QQ^{T}=Q^{T}Q=I_{n}$. Orthogonal matrices satisfy $|Qv|=|v|$ for all
$v\in\Real^{n}$.  Given $C\in\Real^{m\times n}$, $B\in\Real^{n\times
m}$, and $A\in\Real^{n\times n}$, pair $(C,\,A)$ is {\em observable}
if $[C^{T}\, A^{T}C^{T}\, A^{2T}C^{T}\,\ldots\, A^{(n-1)T}C^{T}]$ is
full row rank. Pair $(C,\,A)$ is {\em detectable} (in the
discrete-time sense) if that $CA^{k}x=0$ for some $x\in\Real^{n}$ and
for all $k\in\Natural$ implies $\lim_{k\to\infty}A^{k}x=0$. Pair
$(A,\,B)$ is {\em stabilizable} if $(B^{T},\,A^{T})$ is
detectable. Matrix $A\in\Real^{n\times n}$ is {\em neutrally stable}
(in the discrete-time sense) if it has no eigenvalue with magnitude
greater than unity and the Jordan block corresponding to an eigenvalue
$\lambda$ with $|\lambda|=1$ is of size one.\footnote{Note that $A$ is
neutrally stable iff there exists a symmetric positive definite matrix
$P$ such that $A^{T}PA-P\leq 0$, \cite{antsaklis97}.}  Let
$\one\in\Real^{p}$ denote the vector with all entries equal to unity.

{\em Kronecker product} of $A\in\Real^{m\times n}$ and $B\in\Real^{p\times q}$ is
\begin{eqnarray*}
A\otimes B:=
\left[
\begin{array}{ccc}
a_{11}B & \cdots & a_{1n}B\\
\vdots  & \ddots & \vdots\\
a_{m1}B & \cdots & a_{mn}B
\end{array}
\right]
\end{eqnarray*}
In the pages to come we will enjoy the properties $(A\otimes
B)(C\otimes D)=(AC)\otimes(BD)$ (provided that products $AC$ and $BD$
are allowed), $A\otimes B+A\otimes C=A\otimes(B+C)$ (for $B$ and
$C$ that are of same size) and $|A\otimes B|=|A||B|$.

Matrix $P\in\Real^{n\times n}$ is an {\em orthogonal projection} onto
the subspace ${\rm range}(P)$ if $P^2=P$ and $P^{T}=P$. For an
orthogonal projection $P$, if the columns of $C^{T}\in\Real^{n\times
m}$ are an orthonormal basis for ${\rm range}(P)$ then
$P=C^{T}C$. Matrix $V=I_{n}-P$ is also an orthogonal projection and 
${\rm range}(V)={\rm range}(P)^{\perp}$. It is easy to see that $PV=VP=0$.

A ({\em directed}) {\em graph} is a pair $(\N,\,\A)$ where $\N$ is a
nonempty finite set (of {\em nodes}) and $\A$ is a finite collection
of pairs ({\em arcs}) $(n_{i},\,n_{j})$ with $n_{i},\,n_{j}\in\N$. A
{\em path} from $n_{1}$ to $n_{\ell}$ is a sequence of nodes
$\{n_{1},\,n_{2},\,\ldots,\,n_{\ell}\}$ such that $(n_{i},\,n_{i+1})$
is an arc for $i\in\{1,\,2,\,\ldots,\,\ell-1\}$. A graph is {\em
connected} if it has a node to which there exists a path from every
other node.\footnote{Note that this definition of connectedness for
directed graphs is weaker than strong connectivity and stronger than
weak connectivity.}

The graph of a matrix $\Lambda:=[\lambda_{ij}]\in\Real^{p\times p}$
is the pair $(\N,\,\A)$ where $\N
=\{n_{1},\,n_{2},\,\ldots,\,n_{p}\}$ and $(n_{i},\,n_{j})\in\A$ iff
$\lambda_{ij}>0$. Matrix $\Lambda$ is said to be {\em connected} (in the discrete-time sense) 
if it satisfies:
\begin{enumerate}
\item[(i)] $\lambda_{ii}>0$ and $\lambda_{ij}\geq 0$ for all $i,\,j$;
\item[(ii)] each row sum equals 1;
\item[(iii)] its graph is connected.\footnote{
Recall that for continuous time applications, definition of
connectedness is different: a matrix $[\gamma_{ij}]$ is considered
connected (in the continuous-time sense) if $\gamma_{ij}\geq 0$ for
$i\neq j$; each row sum equals 0; and its graph is connected.  }
\end{enumerate}
For $\Lambda$ that is connected, it is known that
$\lim_{k\to\infty}\Lambda^{k}=\one{r}^{T}$ where $r\in\Real^{p}$ has
nonnegative entries and satisfies $r^{T}\one=1$ and
$r^{T}\Lambda=r^{T}$. We mention that, in an interconnection of
systems, if the matrix describing the network topology satisfies
properties (i) and (ii) above, then the coupling between the systems
is said to be {\em diffusive}.

Given maps
$\xi_{i}:\Natural\to\Real^{n}$ for $i\in\{1,\,2,\,\ldots,\,p\}$ and a map
$\bar\xi:\Natural\to\Real^{n}$, the elements of the set
$\{\xi_{i}(\cdot):i=1,\,2,\,\ldots,\,p\}$ are
said to {\em synchronize to} $\bar{\xi}(\cdot)$ if
$|\xi_{i}(k)-\bar\xi(k)|\to 0$ as $k\to\infty$ for all $i$. 

\section{Problem statement}\label{sec:problemstatement}  

\subsection{Systems under study}
We consider $p$ identical discrete-time linear systems
\begin{eqnarray}\label{eqn:system}
x_{i}^{+}=Ax_{i}+u_{i}\, ,\quad y_{i}=Cx_{i}\, ,\quad i=1,\,2,\,\ldots,\,p
\end{eqnarray}
where $x_{i}\in\Real^{n}$ is the {\em state}, $x_{i}^{+}$ is the state
at the next time instant, $u_{i}\in\Real^{n}$ is the {\em input}, and
$y_{i}\in\Real^{m}$ is the {\em output} of the $i$th system. Matrices
$A$ and $C$ are of proper dimensions. The solution of $i$th system at
time $k\in\Natural$ is denoted $x_{i}(k)$. In this paper we consider
the case where at each time instant only the following information
\begin{eqnarray}\label{eqn:z}
z_{i}&=&\sum_{j=1}^{p}\lambda_{ij}(y_{j}-y_{i})
\end{eqnarray}
is available to $i$th system to determine an input value where
$\Lambda:=[\lambda_{ij}]\in\Real^{p\times p}$ is the matrix describing
the network topology. Matrix $\Lambda$ has nonnegative entries,
strictly positive diagonal entries, and rows summing up to one. That
is, the coupling between systems is diffusive.

\subsection{Assumptions made}
We make the following assumptions on systems~\eqref{eqn:system} which
will henceforth hold.
\\

\noindent
{\bf (A1)} $A$ is neutrally stable.\\
{\bf (A2)} $(C,\,A)$ is detectable.

\subsection{Objectives}
We have two objectives in this paper. The first one is to find the
answer to the question {\em Does there exist a linear feedback law
$L\in\Real^{n\times m}$ such that solutions of
systems~\eqref{eqn:system} with $u_{i}=Lz_{i}$, where $z_{i}$ is as in
\eqref{eqn:z}, globally synchronize to some bounded trajectory 
for all connected $\Lambda$?} Our second objective is, if the answer
to the previous question is affirmative, to {\em devise an algorithm
to compute one such $L$.}

\section{A preliminary result}\label{sec:pre}

Consider $p$ interconnected systems 
\begin{eqnarray}\label{eqn:systemdt}
\xi_{i}^{+}=Q\xi_{i}+QH^{T}H\sum_{j=1}^{p}\lambda_{ij}(\xi_{j}-\xi_{i})\, ,\quad 
i=1,\,2,\,\ldots,\,p
\end{eqnarray}
with $Q\in\Real^{n\times n}$ and
$H\in\Real^{m\times n}$. We make the following assumptions on systems~\eqref{eqn:systemdt} 
which will henceforth hold.\\

\noindent
{\bf (B1)}\ $Q$ is orthogonal.\\
{\bf (B2)}\ $HH^{T}=I_{m}$.\\
{\bf (B3)}\ $(H,\,Q)$ is observable.
\\

\noindent The following result can be
implicitly found, for instance, in \cite[Cor.~15]{tuna07}.
\begin{lemma}\label{lem:nolcos}
Given $\Lambda\in\Real^{p\times p}$ that is connected,
let $r\in\Real^{p}$ be such that
$\lim_{k\to\infty}\Lambda^{k}=\one r^{T}$. Then
there exist $c\geq 1$ and $\sigma\in(0,\,1)$ such that 
$|\Lambda^{k}-\one r^{T}|\leq c\sigma^{k}$ for all $k\in\Natural$.
\end{lemma}

\noindent 
We also need the following result for later use.
\begin{lemma}\label{lem:alpha}
Pair $(H,\,Q)$ of \eqref{eqn:systemdt} satisfies
\begin{eqnarray}\label{eqn:attenuate}
\left|\prod_{i=0}^{n-1}(I_{n}-Q^{iT}H^{T}HQ^{i})\right|<1\,.
\end{eqnarray}
\end{lemma}

\begin{proof}
Note that $Q^{iT}H^{T}HQ^i$ is an orthogonal projection for
$i=0,\,1,\,\ldots,\,n-1$. Whence
$|I_{n}-Q^{iT}H^{T}HQ^i|\in\{0,\,1\}$. Now, suppose
\eqref{eqn:attenuate} is not true.  Then there exists $v\in\Real^{n}$
with $|v|=1$ such that
\begin{eqnarray*}
|(I_{n}-Q^{(n-1)T}H^{T}HQ^{n-1})\cdots(I_{n}-H^{T}H)v|=1\,.
\end{eqnarray*}
For $w\in\Real^{n}$ and an orthogonal projection $P\in\Real^{n\times
n}$, if $(I_{n}-P)w\neq w$ then it must be that $|(I_{n}-P)w|<|w|$. As
a consequence we must have $Q^{iT}H^{T}HQ^iv=0$ for all $i$. Thence $HQ^{i}v=0$ for all $i$. 
This means that $v$ is orthogonal to every column vector of $[H^{T}\ Q^{T}H^{T}\
Q^{2T}H^{T}\ \ldots\ Q^{(n-1)T}H^{T}]$. This however is a
contradiction for $(H,\,Q)$ pair is observable.
\end{proof}\\

\noindent
We now provide a key result.

\begin{theorem}\label{thm:longsought}
Consider systems~\eqref{eqn:systemdt}. Suppose $\Lambda$ is connected
and let $r\in\Real^{p}$ be such that
$\lim_{k\to\infty}\Lambda^{k}=\one r^{T}$. Then solutions
$\xi_{i}(\cdot)$ for $i=1,\,2,\,\ldots,\,p$ synchronize to
\begin{eqnarray*}
\bar{\xi}(k):=({r^{T}\otimes Q^{k}})
\left[\!\!
\begin{array}{c}
\xi_{1}(0)\\
\vdots\\
\xi_{p}(0)
\end{array}
\!\!
\right]
\end{eqnarray*}
\end{theorem}

\begin{proof}
Let us stack individual system states to obtain $\ex:=[\xi_{1}^{T}\ \xi_{2}^{T}\ 
\ldots\ \xi_{p}^{T}]^{T}$.  From \eqref{eqn:systemdt} we obtain
\begin{eqnarray}\label{eqn:fiufiu}
\ex^{+}=(I_{p}\otimes Q)(I_{p}\otimes I_{n}+(\Lambda-I_{p})\otimes(H^{T}H))\ex\,.
\end{eqnarray}
Let $\dabilyu(k):=(I_{p}\otimes Q^{-k})\ex(k)$, $P_{k}:=Q^{kT}H^{T}HQ^{k}$, and
$V_{k}:=I_{n}-P_{k}$ for $k\in\Natural$. Observe that for each $k$, $P_{k}$ is an
orthogonal projection onto subspace ${\rm range}(Q^{kT}H^{T})$
and $V_{k}$ onto ${\rm range}(Q^{kT}H^{T})^{\perp}$.
From \eqref{eqn:fiufiu} we can write
\begin{eqnarray}\label{eqn:unutulmaz}
\dabilyu(k+1)=(I_{p}\otimes V_{k}+\Lambda\otimes P_{k})\dabilyu(k)\,.
\end{eqnarray}
For $k,\,h\in\Natural$ with $k\geq h$, let
\begin{eqnarray*}
\Phi(k,\,h):=\prod_{\tau=h}^{k-1}(I_{p}\otimes V_{\tau}+\Lambda\otimes P_{\tau})
\end{eqnarray*}
with $\Phi(h,\,h)=I_{pn}$. Note that
$\dabilyu(k)=\Phi(k,\,h)\dabilyu(h)$. We now establish the following:
\begin{eqnarray}\label{eqn:claim}
\lim_{k\to\infty}\Phi(k,\,h)=\one r^{T}\otimes I_{n}
\end{eqnarray}
for any fixed $h$. Without loss of generality let $h=0$.
For $\ell,\,k\in\Natural$ with $\ell\leq k$ let us define
\begin{eqnarray}\label{eqn:defineMlk}
M_{\ell,k}:=\sum_{M\in\Omega_{\ell,k}}M
\end{eqnarray}
where
\begin{eqnarray*}
\Omega_{\ell,k}:=\{M:M=L_{k-1}L_{k-2}\cdots L_{0},
L_{i}\in\{V_{i},\,P_{i}\},\,\#\{i:L_{i}=P_{i}\}=\ell\}\,.
\end{eqnarray*} 
For instance, $\Omega_{0,4}=\{V_{3}V_{2}V_{1}V_{0}\}$ and 
$\Omega_{2,4}=\{V_{3}V_{2}P_{1}P_{0},$ $\,V_{3}P_{2}V_{1}P_{0},\,
P_{3}V_{2}V_{1}P_{0},\,$
$V_{3}P_{2}P_{1}V_{0},\,P_{3}V_{2}P_{1}V_{0},\,P_{3}P_{2}V_{1}V_{0}\}$. Observe
that
\begin{eqnarray}\label{eqn:omegacard}
\#\Omega_{\ell,k}=\frac{k!}{(k-\ell)!\ell!}=:\left(\!\begin{array}{c}k\\ \ell\end{array}\!\right)
\end{eqnarray}
Note that $M_{\ell,k+1}=V_{k}M_{\ell,k}+P_{k}M_{\ell-1,k}$ for $\ell\in\{1,\,2,\,\ldots,\,k\}$, 
$M_{0,k+1}=V_{k}M_{0,k}$, and $M_{k+1,k+1}=P_{k}M_{k,k}$. We can write
\begin{eqnarray}\label{eqn:whatphiis}
\Phi(k,\,0)=\sum_{\ell=0}^{k}(\Lambda^{\ell}\otimes M_{\ell,k})
\end{eqnarray}
and 
\begin{eqnarray}\label{eqn:Msum}
\sum_{\ell=0}^{k}M_{\ell,k}=I_{n}\,.
\end{eqnarray}
Let $\alpha:=|V_{n-1}V_{n-2}\cdots V_{0}|$. Lemma~\ref{lem:alpha} guarantees that $\alpha<1$.
We make the following observations. For all $\ell,\,k$, 
\begin{eqnarray}\label{eqn:observe1}
|M_{\ell,k}|\leq 1
\end{eqnarray}
and
\begin{eqnarray}\label{eqn:observe2}
|M_{\ell,k}|\leq \left(\!\begin{array}{c}k\\ \ell\end{array}\!\right)
\alpha^{\lfloor\frac{k+1}{n}\rfloor-\ell}
\end{eqnarray} 
for $\lfloor (k+1)/n\rfloor-\ell\geq 0$.  

Let us first show
\eqref{eqn:observe1}. Suppose for some $k\in\Natural$ and all
$v\in\Real^{n}$ with $|v|=1$ we have
\begin{eqnarray}\label{eqn:somek}
\sum_{\ell=0}^{k}|M_{\ell,k}v|^{2}=1\,.
\end{eqnarray}
Then we can write, since both $V_{k}$ and $P_{k}$ are orthogonal
projections satisfying $V_{k}P_{k}=0$,
\begin{eqnarray*}
\sum_{\ell=0}^{k+1}|M_{\ell,k+1}v|^{2}
&=& |V_{k}M_{0,k}v|^{2}+|P_{k}M_{k,k}v|^{2}
+\sum_{\ell=1}^{k}|V_{k}M_{\ell,k}v+P_{k}M_{\ell-1,k}v|^{2}\\
&=& |V_{k}M_{0,k}v|^{2}+|P_{k}M_{k,k}v|^{2}
+\sum_{\ell=1}^{k}\left(|V_{k}M_{\ell,k}v|^{2}+|P_{k}M_{\ell-1,k}v|^{2}\right)\\
&=&\sum_{\ell=0}^{k}\left(|V_{k}M_{\ell,k}v|^{2}+|P_{k}M_{\ell,k}v|^{2}\right)\\
&=&\sum_{\ell=0}^{k}|M_{\ell,k}v|^{2}\\
&=& 1\,.
\end{eqnarray*}
Hence, by induction, \eqref{eqn:somek} holds for all $k$ since it trivially holds for $k=0$ thanks to 
\eqref{eqn:Msum}. We therefore have \eqref{eqn:observe1} as a direct implication of \eqref{eqn:somek}.

Now, we show \eqref{eqn:observe2}. Let us be given some $k\geq
n-1$. Let then $k_{0}:= k-n+1$. We observe that
\begin{eqnarray*}
|V_{k}V_{k-1}\cdots V_{k-n+1}|
&=& |Q^{k_{0}T}V_{n-1}Q^{k_{0}}Q^{k_{0}T}V_{n-2}Q^{k_{0}}\cdots Q^{k_{0}T}V_{0}Q^{k_{0}}|\\
&=& |Q^{k_{0}T}V_{n-1}V_{n-2}\cdots V_{0}Q^{k_{0}}|\\
&=& |V_{n-1}V_{n-2}\cdots V_{0}|\\
&=& \alpha\,.
\end{eqnarray*}
Hence for $M\in\Omega_{\ell,k}$ one can write
\begin{eqnarray*}
|M|\leq \alpha^{\lfloor\frac{k+1}{n}\rfloor-\ell}
\end{eqnarray*}
for $\lfloor(k+1)/n\rfloor-\ell\geq 0$. Then by \eqref{eqn:defineMlk} and \eqref{eqn:omegacard} 
we obtain \eqref{eqn:observe2}.

Now we are ready to show \eqref{eqn:claim}. Let us be given
$\delta>0$. Let $c\geq 1$ and $\sigma\in(0,\,1)$ be such that $|\Lambda^{k}-\one r^{T}|\leq c\sigma^{k}$ for all $k\in\Natural$. 
Such $c$ and $\sigma$ exist by Lemma~\ref{lem:nolcos}. Choose $\ell^{*}\in\Natural$ such that
$\sum_{\ell=\ell^{*}}^{\infty}c\sigma^{\ell}\leq\delta/2$. Then choose
$k^{*}\in\Natural$ such that
\begin{eqnarray*}
\sum_{\ell=0}^{\ell^{*}-1}\left(\!\begin{array}{c}k\\ \ell\end{array}\!\right)
\alpha^{\lfloor\frac{k+1}{n}\rfloor-\ell}\leq\frac{\delta}{2c}
\end{eqnarray*}
for all $k\geq k^{*}$. Now, let us be given some $k\geq k^{*}$. We
write by \eqref{eqn:whatphiis} and \eqref{eqn:Msum}
\begin{eqnarray*}
|\Phi(k,\,0)-\one r^{T}\otimes I_{n}|
&=& \left|\sum_{\ell=0}^{k}\Lambda^{\ell}\otimes M_{\ell,k}-\sum_{\ell=0}^{k}\one r^{T}
\otimes M_{\ell,k}\right|\\
&=& \left|\sum_{\ell=0}^{k}(\Lambda^{\ell}-\one r^{T})\otimes M_{\ell,k}\right|\\
&\leq& \sum_{\ell=0}^{\ell^{*}-1}\left|(\Lambda^{\ell}-\one r^{T})\otimes M_{\ell,k}\right|
+\sum_{\ell=\ell^{*}}^{k}\left|(\Lambda^{\ell}-\one r^{T})\otimes M_{\ell,k}\right|\\
&=& \sum_{\ell=0}^{\ell^{*}-1}|\Lambda^{\ell}-\one r^{T}||M_{\ell,k}|
+\sum_{\ell=\ell^{*}}^{k}|\Lambda^{\ell}-\one r^{T}||M_{\ell,k}|\\
&\leq& c\sum_{\ell=0}^{\ell^{*}-1}|M_{\ell,k}|
+\sum_{\ell=\ell^{*}}^{k}|\Lambda^{\ell}-\one r^{T}|\\
&\leq& c\sum_{\ell=0}^{\ell^{*}-1}\left(\!\begin{array}{c}k\\ \ell\end{array}\!\right)
\alpha^{\lfloor\frac{k+1}{n}\rfloor-\ell}
+\sum_{\ell=\ell^{*}}^{k}c\sigma^{\ell}\\
&\leq& c\frac{\delta}{2c}
+\frac{\delta}{2}\\
&=&\delta\
\end{eqnarray*}
where we have employed \eqref{eqn:observe1} and \eqref{eqn:observe2}. 
Having shown \eqref{eqn:claim}, we can write
\begin{eqnarray}\label{eqn:JP}
\lim_{k\to\infty}|\ex(k)-(\one r^{T}\otimes Q^{k})\ex(0)|
&=&\lim_{k\to\infty}|(I_{p}\otimes Q^{k})\dabilyu(k)-(\one r^{T}\otimes Q^{k})\dabilyu(0)|\nonumber\\
&=&\lim_{k\to\infty}|(I_{p}\otimes Q^{k})\Phi(k,\,0)\dabilyu(0)-(\one r^{T}\otimes Q^{k})\dabilyu(0)|\nonumber\\
&=&\lim_{k\to\infty}|((I_{p}\otimes Q^{k})(\one r^{T}\otimes I_{n})-\one r^{T}\otimes Q^{k})\dabilyu(0)|\nonumber\\
&=&\lim_{k\to\infty}|(\one r^{T}\otimes Q^{k}-\one r^{T}\otimes Q^{k})\dabilyu(0)|\nonumber\\
&=&0
\end{eqnarray}
where we used the fact that $\ex(0)=\dabilyu(0)$. Eq.~\eqref{eqn:JP} implies
\begin{eqnarray*}
\lim_{k\to\infty}|\xi_{i}(k)-(r^{T}\otimes Q^{k})\ex(0)|=0
\end{eqnarray*}
for all $i\in\{1,\,2,\,\ldots,\,p\}$. Hence the result.
\end{proof}

\section{Synchronization via output feedback}\label{sec:main}

We are now ready to answer the question asked in the problem
statement: Does there exist $L\in\Real^{m\times n}$ such that
solutions of systems~\eqref{eqn:system} with $u_{i}=Lz_{i}$, where
$z_{i}$ is as in \eqref{eqn:z}, synchronize for all connected
$\Lambda$? The answer, we will see, is affirmative and lies in a
straightforward generalization of the key result
(Theorem~\ref{thm:longsought}) of the previous section. We also
provide a simple algorithm to calculate such $L$. Let us begin with
the following fact.

\begin{fact}\label{fact:one}
Let $F\in\Real^{n\times n}$ be a neutrally stable matrix with all its eigenvalues 
having unity magnitude. Then there exists $\Real^{n\times n}\ni R=R^{T}>0$ such that $F^{T}RF=R$.
\end{fact}

\begin{proof}
Since $F$ has no Jordan block of size greater than one, it can be
diagonalized. Therefore there exist $Z\in\Complex^{n\times n}$ and a
diagonal matrix $D\in\Complex^{n\times n}$ such that
$F=ZDZ^{-1}$. Since the diagonal entries of $D$ are all of unity
magnitude, $|D^{k}x|=|x|$ for all $x\in\Real^{n}$ and $k=1,\,2,\,\ldots$.
Therefore there exist real numbers $0<a\leq b<\infty$ such that 
\begin{eqnarray*}
a|x|\leq|F^{k}x|\leq b|x|\,,\qquad k=1,\,2,\,\ldots
\end{eqnarray*}
Note that $a^{2}x^{T}x\leq x^{T}F^{kT}F^{k}x\leq b^{2}x^{T}x$ for all
$x$ and $k$. Let us define the compact set $\X:=\{X\in\Real^{n\times
n}: X=X^{T}, a^{2}I_{n}\leq X\leq b^{2}I_{n}\}$ and the continuous
function $f:\X\to\Real$ as $f(X):=|F^{T}XF-X|$. Finally let
\begin{eqnarray*}
X_{k}:=k^{-1}\sum_{i=1}^{k}F^{iT}F^{i}\,,\qquad k=1,\,2,\,\ldots
\end{eqnarray*} 
By construction $a^{2}I_{n}\leq X_{k}\leq b^{2}I_{n}$ and $X_{k}^{T}=X_{k}$ for all
$k$. Hence $X_{k}\in\X$. We now can write
\begin{eqnarray*}
f(X_{k})
&=&|F^{T}X_{k}F-X_{k}|\\
&=&k^{-1}\left|\sum_{i=2}^{k+1}F^{iT}F^{i}-\sum_{i=1}^{k}F^{iT}F^{i}\right|\\
&=&k^{-1}|F^{(k+1)T}F^{k+1}-F^{T}F|\\
&\leq&k^{-1}2b^{2}\,.
\end{eqnarray*}
As a result we have $\lim_{k\to\infty}f(X_{k})=0$. Since $f$ takes only nonnegative values we deduce
\begin{eqnarray*}
\inf_{X\in\X} f(X)=0\,.
\end{eqnarray*}
Compactness of $\X$ together with continuity of $f$ implies that minimum is attained 
\cite[Cor.~6.57]{browder96}. Thus there exists 
$R\in\X$ such that $f(R)=0$. 
\end{proof}

\begin{algorithm}\label{alg:L}
Given $A\in\Real^{n\times n}$ that is neutrally stable and $C\in\Real^{m\times n}$,   
we obtain $L\in\Real^{n\times m}$ as follows. Let $n_{1}\leq n$ be the number of 
eigenvalues of $A$ with unity magnitude. Let $n_{2}:=n-n_{1}$. If $n_{1}=0$, then let $L:=0$; else
construct $L$ according to the following steps. \\

\noindent
{\em Step 1:} Choose $U\in\Real^{n\times n_{1}}$ and
$W\in\Real^{n\times n_{2}}$ satisfying
\begin{eqnarray*}
[U\ W]^{-1}A[U\ W]=
\left[
\begin{array}{cc}
F & 0\\
0 & G
\end{array}
\right]
\end{eqnarray*} 
where all the eigenvalues of $F\in\Real^{n_{1}\times n_{1}}$ have
unity magnitude. (Assume, without loss of generality for our purposes,
that $CU$ is full row rank.)\\

\noindent
{\em Step 2:} Choose $R\in\Real^{n_{1}\times n_{1}}$ with $R=R^{T}>0$ such that
$F^{T}RF=R$. (This we can do thanks to Fact~\ref{fact:one}.)\\

\noindent
{\em Step 3:} Choose $H\in\Real^{m\times n_{1}}$ satisfying ${\rm
range}(H^{T})={\rm range}(R^{-1/2}U^{T}C^{T})$ and
$HH^{T}=I_{m}$. (Note that then matrix $CUR^{-1/2}H^{T}$ is
invertible.)\\

\noindent
{\em Step 4:} Define $L:=UFR^{-1/2}H^{T}(CUR^{-1/2}H^{T})^{-1}$.
\end{algorithm} 

\noindent
Below is our main result.
\begin{theorem}\label{thm:main}
Consider systems \eqref{eqn:system}. Let $u_{i}=Lz_{i}$ where
$L\in\Real^{n\times m}$ is constructed according to
Algorithm~\ref{alg:L} and $z_{i}$ is as in \eqref{eqn:z}. Then for all
network topologies described by connected $\Lambda$, solutions
$x_{i}(\cdot)$, for $i=1,\,2,\,\ldots,\,p$, synchronize to 
\begin{eqnarray*}
\bar{x}(k):=(r^{T}\otimes A^{k})
\left[
\begin{array}{c}
x_{1}(0)\\
\vdots\\
x_{p}(0)
\end{array}
\right]
\end{eqnarray*}
where $r\in\Real^{p}$ is such that $r^{T}\Gamma=r^{T}$ and $r^{T}\one=1$.
\end{theorem}

\begin{proof}
Let the variables that are not introduced here be defined as in
Algorithm~\ref{alg:L}.  Without loss of generality we assume that $CU$
is full row rank.  Since $H^{T}H$ is an orthogonal projection onto
${\rm range}(H^{T})= {\rm range}(R^{-1/2}U^{T}C^{T})$, we can write
$H^{T}HR^{-1/2}U^{T}C^{T}=R^{-1/2}U^{T}C^{T}$.  Taking the transpose
we obtain $CUR^{-1/2}H^{T}H=CUR^{-1/2}$.  Since $CUR^{-1/2}H^{T}$ is
invertible we obtain $H=(CUR^{-1/2}H^{T})^{-1}CUR^{-1/2}$.  Therefore
$LCUR^{-1/2}=UFR^{-1/2}H^{T}H$. Also, detectability of $(C,\,A)$
implies that pair $(H,\,Q)$ is observable for $Q:=R^{1/2}FR^{-1/2}$.
Note that $Q$ is orthogonal due to $F^{T}RF=R$.

We let $U^{\dagger}\in\Real^{n_{1}\times n}$ and
$W^{\dagger}\in\Real^{n_{2}\times n}$ be such that
\begin{eqnarray*}
\left[
\begin{array}{c}
U^{\dagger}\\
W^{\dagger}
\end{array}
\right]=[U\ W]^{-1}\,.
\end{eqnarray*} 
Note then that $U^{\dagger}U=I_{n_{1}}$, $W^{\dagger}W=I_{n_{2}}$,
$U^{\dagger}W=0$, and $W^{\dagger}U=0$.  Since $u_{i}=Lz_{i}$, we can
combine \eqref{eqn:system} and \eqref{eqn:z} to obtain
$x_{i}^{+}=Ax_{i}+LC\sum_{j=1}^{p}\lambda_{ij}(x_{j}-x_{i})$. By
change of variables $\xi_{i}:=[R^{1/2}\ 0][U\ W]^{-1}x_{i}$ and
$\eta_{i}:=[0\ I_{n_{2}}][U\ W]^{-1}x_{i}$ we can write
\begin{eqnarray}
\xi_{i}^{+}\!\!&=&\!\!Q\xi_{i}
+QH^{T}H\sum_{j=1}^{p}\lambda_{ij}(\xi_{j}-\xi_{i})
+R^{1/2}U^{\dagger}LCW\sum_{j=1}^{p}\lambda_{ij}(\eta_{j}-\eta_{i})\label{eqn:badem1}\\
\eta_{i}^{+}\!\!&=&\!\!G\eta_{i}\,.\label{eqn:badem2}
\end{eqnarray}

Let $\Lambda$ be connected and $r\in\Real^{p}$ be such that
$\lim_{k\to\infty}\Lambda^{k}=\one{r^{T}}$. Then define $\omega_{i}:\Natural\to\Real^{n_{1}}$ 
as $\omega_{i}(k):=Q^{-k}\xi_{i}(k)$ and
$i=1,\,2,\,\ldots,\,p$. Let $\dabilyu:=[\omega_{1}^{T}\ \omega_{2}^{T}\ \ldots\
\omega_{p}^{T}]^{T}$ and $\vi:=[\eta_{1}^{T}\ \eta_{2}^{T}\ \ldots\
\eta_{p}^{T}]^{T}$. Starting from \eqref{eqn:badem1} and
\eqref{eqn:badem2} we can write
\begin{eqnarray*}
\dabilyu(k+1)&=&\left(I_{pn_{1}}+(\Lambda-I_{p})\otimes Q^{-k}H^{T}HQ^{k})\right)\dabilyu(k)\nonumber\\
&&\qquad+\left((\Lambda-I_{p})\otimes Q^{-k-1}MG^{k}\right)\vi(0)
\end{eqnarray*}
where $M:=R^{1/2}U^{\dagger}LCW$.
Thence
\begin{eqnarray}\label{eqn:integral}
\dabilyu(k)=\Phi(k,\,0)\dabilyu(0)
+\left[\sum_{\ell=0}^{k-1}\Phi(k,\,\ell+1)\left((\Lambda-I_{p})\otimes Q^{-\ell-1}MG^{\ell}
\right)
\right]\vi(0)
\end{eqnarray}
where 
\begin{eqnarray*}
\Phi(k,\,\ell):=\prod_{\tau=\ell}^{k-1}\left(I_{pn_{1}}+(\Lambda-I_{p})\otimes Q^{-\tau}H^{T}HQ^{\tau}
\right)
\end{eqnarray*}
is the state transition matrix \cite{antsaklis97}. From
Theorem~\ref{thm:longsought} we can deduce that $\Phi(k,\,\ell)$ is
uniformly bounded for all $k$ and $\ell$. Also, for any fixed $\ell$
we have $\lim_{k\to\infty}\Phi(k,\,\ell)=\one{r^{T}}\otimes
I_{n_{1}}$. Moreover, $Q^{k}$ is uniformly bounded for all $k$, and
$G^{k}$ decays exponentially as $k\to\infty$ for all the eigenvalues
of $G$ are strictly within the unit circle. Therefore we can write
\begin{eqnarray*}
\lefteqn{
\lim_{k\to\infty}\sum_{\ell=0}^{k-1}\Phi(k,\,\ell+1)\left((\Lambda-I_{p})\otimes Q^{-\ell-1}MG^{\ell}
\right)}\\
&&=\sum_{\ell=0}^{\infty}\left(\lim_{k\to\infty}\Phi(k,\,\ell+1)\right)
\left((\Lambda-I_{p})\otimes Q^{-\ell-1}MG^{\ell}
\right)\\
&&=\sum_{\ell=0}^{\infty}(\one{r^{T}}\otimes I_{n_{1}})
\left((\Lambda-I_{p})\otimes Q^{-\ell-1}MG^{\ell}
\right)\\
&&=0\,.
\end{eqnarray*}
Then, by \eqref{eqn:integral}, we can write
\begin{eqnarray*}
\lim_{k\to\infty}\dabilyu(k)=(\one{r^{T}}\otimes I_{n_{1}})\dabilyu(0)\,.
\end{eqnarray*}
Therefore solutions $\xi_{i}(\cdot)$ synchronize to $(r^{T}\otimes
Q^{k})\dabilyu(0)$. Moreover, by \eqref{eqn:badem2}, there is no harm
in claiming that solutions $\eta_{i}(\cdot)$ synchronize to
$(r^{T}\otimes G^{k})\vi(0)$ for $G^k$ is decaying as
$k\to\infty$. We can conclude that solutions $x_{i}(\cdot)$
synchronize to
\begin{eqnarray*}
\lefteqn{
\left(r^{T}\otimes 
\left[UR^{-1/2}\ \ W\right]
\left[
\begin{array}{cc}
Q^{k}&0\\0&G^{k}
\end{array}
\right]
\left[\begin{array}{c}R^{1/2}U^{\dagger}\\ W^{\dagger}\end{array}\right]
\right)
\left[
\begin{array}{c}
x_{1}(0)\\
\vdots\\
x_{p}(0)
\end{array}
\right]}\\
&&\hspace{3in}=(r^{T}\otimes A^{k})
\left[
\begin{array}{c}
x_{1}(0)\\
\vdots\\
x_{p}(0)
\end{array}
\right]
\end{eqnarray*}
Hence the result.
\end{proof}

\section{Dual problem}\label{sec:dual}
In this section we present a problem similar, in fact dual, to the one
stated in Section~\ref{sec:problemstatement}. Consider $p$ identical systems
\begin{eqnarray}\label{eqn:systemdual}
x_{i}^{+}=A^{T}x_{i}+C^{T}u_{i}\, ,\quad i=1,\,2,\,\ldots,\,p
\end{eqnarray}
where $x_{i}\in\Real^{n}$ is the state and $u_{i}\in\Real^{m}$ is the
input of the $i$th system. Matrices $A^{T}$ and $C^{T}$ are of proper
dimensions. Let pair $(A^{T},\,C^{T})$ be stabilizable. Suppose now
that at each time instant the following information
\begin{eqnarray}\label{eqn:zdual}
z_{i}&=&\sum_{j=1}^{p}\lambda_{ij}(x_{j}-x_{i})
\end{eqnarray}
is available to $i$th system to determine an input value. Now the obvious
question we ask is the following. {\em Can we design a linear feedback law 
$K\in\Real^{m\times n}$ such that solutions of
systems~\eqref{eqn:systemdual} with $u_{i}=Kz_{i}$, where $z_{i}$ is as in
\eqref{eqn:zdual}, globally synchronize to some bounded trajectory
for all connected $\Lambda$?} The answer is the next result which 
follows from Theorem~\ref{thm:main}.

\begin{theorem}\label{thm:maindual}
Consider systems \eqref{eqn:systemdual}. Let $u_{i}=L^{T}z_{i}$ where
$L\in\Real^{n\times m}$ is constructed according to
Algorithm~\ref{alg:L} and $z_{i}$ is as in \eqref{eqn:zdual}. Then for all
network topologies described by connected $\Lambda$, solutions
$x_{i}(\cdot)$, for $i=1,\,2,\,\ldots,\,p$, synchronize to 
\begin{eqnarray*}
\bar{x}(k):=(r^{T}\otimes A^{kT})
\left[
\begin{array}{c}
x_{1}(0)\\
\vdots\\
x_{p}(0)
\end{array}
\right]
\end{eqnarray*}
where $r\in\Real^{p}$ is such that $r^{T}\Gamma=r^{T}$ and $r^{T}\one=1$.
\end{theorem}

\bibliographystyle{plain}         
\bibliography{references}            

\end{document}